\documentclass[11pt]{article} \usepackage{amsmath,amsfonts,latexsym,amssymb,amsthm} \usepackage{url}

 \newcommand{\Q}{\mathbb{Q}}  \newcommand{\Z}{\mathbb{Z}}  \newcommand{\F}{\mathbb{F}}

\begin{document}

\title{Torsion groups of elliptic curves over quadratic fields} \author{Sheldon Kamienny and Filip Najman} \date{} \maketitle \begin{abstract} We describe methods to determine all the possible torsion groups of an elliptic curve that actually appear over a fixed quadratic field. We use these methods to find, for each group that can appear over a quadratic field, the field with the smallest absolute value of its discriminant such that there exists an elliptic curve with that torsion. We also examine the interplay of the torsion and rank over a fixed quadratic field and see that what happens is very different than over $\Q$. Finally we give some results concerning the number and density of fields with an elliptic curve with given torsion over them. \end{abstract} \textbf{Keywords} Torsion group, Elliptic Curves, Quadratic fields \section{Introduction.}

For an elliptic curve $E$ over a number field $K$, it is well known, by the Mordell-Weil theorem, that the set $E(K)$ of $K$-rational points on $E$ is a finitely generated abelian group. The group $E(K)$ is isomorphic to $T\oplus\mathbb \Z^r$, where $r$ is a non-negative integer and $T$ is the torsion subgroup. When $K=\mathbb Q$, by Mazur's Theorem \cite{Maz}, the torsion subgroup is either cyclic of order $m$, where $1 \leq m \leq 10$ or $m=12$, or of the form $\mathbb Z/2 \Z \oplus \mathbb Z/2m \Z$, where $1 \leq m \leq 4$.

If $K$ is a quadratic field, then the following theorem of Kenku and Momose \cite{km} and the fist author \cite{Kam1} classifies the possible torsions. \newtheorem{tm}{Theorem} \begin{tm} Let $K$ be a quadratic field and $E$ an elliptic curve over $K$. Then the torsion subgroup $E(K)_{tors}$ of $E(K)$ is isomorphic to one of the following $26$ groups: $$\mathbb Z/m \Z, \text{ for } 1 \leq m\leq 18,\ m\neq 17,$$ $$\mathbb Z/2 \Z \oplus \mathbb Z/2m \Z, \text{ for } 1 \leq m\leq 6,$$ $$\mathbb Z/3 \Z \oplus \mathbb Z/3m \Z, \text{ for }  m=1,2$$ $$\mathbb Z/4 \Z \oplus \mathbb Z/4 \Z.$$ \end{tm} While this theorem settles the question which torsion groups appear if the field varies through all quadratic fields, it tells us nothing about the possible torsion subgroups if we fix a certain quadratic field.

The second author found all the possible torsions over each of the two cyclotomic quadratic fields, $\Q (i)$ and $\Q (\sqrt{-3})$, in \cite{Naj1} and \cite{Naj2}. The torsion groups appearing over $\Q (i)$ are the ones from Mazur's theorem and $\Z/4 \Z \oplus \Z/4 \Z$, while the torsion groups appearing over $\Q (\sqrt{-3})$ are the ones from Mazur's theorem, $\Z/3 \Z \oplus \Z/3 \Z$ and $\Z/3 \Z \oplus \Z/6 \Z$.

In this paper we describe methods that can be used to obtain results of this type, i.e. find all possible torsions over a given quadratic field. The problem amounts to finding whether certain modular curves that are either elliptic or hyperelliptic have $K$-rational points that are not cusps for the given quadratic field $K$.

We use these techniques to find, for each possible torsion group, the exact field with the smallest absolute value of the discriminant such that an elliptic curve exists with that torsion group. This is done by searching through fields with ascending discriminant until we find one over which the torsion group is possible.

Over the rationals, there are infinitely many nonisomorphic elliptic curves with each torsion group. This is not generally the case for all torsion groups over quadratic fields. Some torsion groups will appear for only finitely many elliptic curves (up to isomorphism) for each field, some will always appear infinitely many times if they appear at all, while some will appear finitely many times over some fields and infinitely many times over others.

A folklore conjecture is that the rank of an elliptic curve over the rationals with any possible torsion can be arbitrarily large. In contrast, we will find the maximum rank that an elliptic curve with prescribed torsion can have over certain quadratic fields. This is done with torsion groups that appear finitely often over the given fields.

\section{Finding the possible torsions over a fixed quadratic field.} Let $K$ be a quadratic field. Denote by $Y_1(m,n)$ the affine curve whose $K$-rational points classify isomorphism classes of the triples $(E, P_m, P_n)$, where $E$ is an elliptic curve (over $K$) and $P_m$ and $P_n$ are torsion points (over $K$) which generate a subgroup isomorphic to $\Z/m \Z \oplus \Z/n \Z$. For simplicity,  we write $Y_1(n)$ instead of $Y_1(1,n)$. Let $X_1(m,n)$ be the compactification of the curve $Y_1(m,n)$ obtained by adjoining its cusps.

Nice models of the curves $X_1(n)$ can be found, for example, in \cite{Baa}, while the curves $X_1(2,n)$ that we will need can be found in \cite{Rab}. We list the curves (and their cusps) that correspond to the torsion points that appear over quadratic field, but not over the rationals. We also exclude the curves $X_1(3,3),\ X_1(3,6)$ and $X_1(4,4)$ since they appear only over the cyclotomic quadratic fields, that are dealt with in \cite{Naj1} and \cite{Naj2}. The curves are as follows: \begin{equation} X_1(11): y^2-y=x^3-x^2, \end{equation} where the cusps satisfy \begin{equation} x(x-1)(x^5-18x^4+35x^3-16x^2-2x+1)=0, \end{equation}

\begin{equation} X_1(13): y^2=x^6-2x^5+x^4-2x^3+6x^2-4x+1, \end{equation} where the cusps satisfy \begin{equation} x(x-1)(x^3-4x^2+x+1)=0, \end{equation}

\begin{equation} X_1(14): y^2+xy+y=x^3-x \end{equation} where the cusps satisfy \begin{equation} x(x-1)(x+1)(x^3-9x^2-x+1)(x^3-2x^2-x+1)=0, \end{equation}

\begin{equation} X_1(15): y^2+xy+y=x^3+x^2 \end{equation} where the cusps satisfy \begin{equation} x(x+1)(x^2+x+1)(x^4+3x^3+4x^2+2x+1)(x^4-7x^3-6x^2+2x+1)=0, \end{equation}

\begin{equation} X_1(16): y^2=x(x^2+1)(x^2+2x-1) \end{equation} where the cusps satisfy \begin{equation} x(x-1)(x+1)(x^2-2x-1)(x^2+2x-1)=0, \end{equation}

\begin{equation} X_1(18): y^2=x^6+2x^5+5x^4+10x^3+10x^2+4x+1 \end{equation} where the cusps satisfy \begin{equation} x(x+1)(x^2+x+1)(x^3-3x-1)=0, \end{equation}

\begin{equation} X_1(2,10): y^2=x^3+x^2-x \end{equation} where the cusps satisfy \begin{equation} x(x-1)(x+1)(x^2+x-1)(x^2-4x-1)=0, \end{equation}

\begin{equation} X_1(2,12): y^2=x^3-x^2+x \end{equation} where the cusps satisfy \begin{equation} x(x-1)(x+1)(x^2+1)(x^2-x+1)(x^2-4x+1)=0. \end{equation}

In order to find whether there exist a curve with torsion $\Z/m \Z \oplus \Z/n \Z$ over a quadratic field $K$, one needs to determine whether $X_1(m,n)$ has a $K$-rational point that is not a cusp.

If $X_1(m,n)$ is an elliptic curve, then a usual method of computing the rank is to perform $2$-descent. One can use an implementation \cite{sim} of Simon in PARI/GP. If the rank is positive then there will be infinitely many elliptic curves with torsion $\Z/m \Z \oplus \Z/n \Z$ over $K$. If the rank is zero, then one has to check whether all the torsion points are cusps. If not, then there will be finitely, explicitly computable elliptic curves with the given torsion subgroup.

If $X_1(m,n)$ is a hyperelliptic curve, there are, by Faltings' theorem, finitely many $K$-rational points, implying that there are finitely many elliptic curves (up to isomorphism) with torsion $\Z/m \Z \oplus \Z/n \Z$ over $K$. To find all the points one can sometimes proceed to compute the rank of the Jacobian using $2$-descent on Jacobians. This can be done in MAGMA (see \cite{st}). Although this will often work, $2$-descent is not an algorithm, as one has no guarantee that one will obtain the exact rank using it, only upper bounds.

Because of this we give an alternative approach, using the method of Mazur in Section 4, that will give us a criterion when the Jacobian of $X_1(13)$ and $X_1(18)$ will have rank $0$.

If the rank is equal to zero, after finding the torsion of the Jacobian, one has to check whether any of the torsion points arise from a $K$-rational point that is not a cusp. If the rank is positive, this significantly complicates the problem and one can try to apply the Chabauty method \cite{sik} (if the rank is 1) or some other similar method.

Note that the only other hyperelliptic curve, $X_1(16)$, is generally easier to deal with, since $f(x)$ (where $y^2=f(x)$ is a model of $X_1(16)$) is not irreducible. This enables one to try to find all the points with more elementary methods, like covering the hyperelliptic curve with 2 elliptic curves (this is essentially what is done in \cite{Naj1}).

Note that once a $K$-rational point on $Y_1(m,n)$ is found, one can find in \cite{Rab} how to actually construct an elliptic curve with torsion $\Z/m \Z \oplus \Z/n \Z$ over $K$.

\section{Smallest field with a given torsion}

We now apply methods described in section 2 to find the field with the smallest $|\Delta |$ over which each torsion group appears, where $\Delta$ is the discriminant. We start with the group $\Z/11 \Z$.

\newtheorem{tm2}[tm]{Theorem} \begin{tm2} \label{t2} The quadratic field $K$ with smallest $|\Delta |$ such that $\Z/11 \Z$ appears as a torsion group over $K$ is $K=\Q (\sqrt{-7})$. \end{tm2} \begin{proof} The torsion group of $X_1(11)(K)$ is isomorphic to $\Z/5 \Z$ for all quadratic fields $K$ (see \cite[Lemme 2.1, p. 26]{Rab}). All the torsion points have $x=0$ or $x=1$, implying that the torsion points are cusps. We now compute, using 2-descent, that $rank(X_1(11)(\Q(\sqrt d)))=0$ for $d=-1,\ -3$ and $5$, and $rank(X_1(11)(\Q(\sqrt {-7}))=1$. Taking a nontorsion point on $X_1(11)(\Q(\sqrt {-7})$, we get the curve $$y^2 + \frac{ 85+33\sqrt{-7}}{128}xy + \frac{85\sqrt{-7}-999}{16384}y = x^3 + \frac{89\sqrt{-7} - 275}{4096}x^2,$$ where $(0,0)$ is a torsion point of order $11$. \end{proof}

Next we examine the group $\Z/13 \Z$.

\newtheorem{tm3}[tm]{Theorem} \begin{tm3} \label{t3} The quadratic field $K$ with smallest $|\Delta |$ such that $\Z/13 \Z$ appears as a torsion group over $K$ is $K=\Q (\sqrt{17})$. \end{tm3} \begin{proof} First note that from \cite{Naj2} we can see that $\Z/13 \Z$ does not appear over $\Q(i)$ and $\Q(\sqrt{-3})$, the two fields with the smallest $|\Delta |$. As $X_1(13)$ is a hyperelliptic curve, we are led to the study of its Jacobian $J_1(13)$. We obtain via 2-descent that $rank(J_1(13)(\Q(\sqrt d)))=0$ for $d=5,-7,2,-2,-11,3,13$ and $-15$. The torsion of $J_1(13)(\Q)$ is isomorphic to $\Z/21 \Z$ and all the points on $J_1(13)(\Q)_{tors}$ are generated by the cusps of $X_1(13)$. We will prove $J_1(13)(\Q(\sqrt d))_{tors}=J_1(13)(\Q)_{tors}$ for $d=5,-7,2,-2,-11,3,13$ and thus complete our proof. Note that if $p$ is a inert (in $\Q(\sqrt d)$) prime of good reduction, then the prime-to-$p$ part of $J_1(13)(\Q(\sqrt d))_{tors}$ injects into $J_1(13)(\F_p(\sqrt d)) \simeq J_1(13)(\F_{p^2})$. If $p$ splits, let $\mathfrak p$ be a prime over $p$. Then $O_K/\mathfrak p \simeq \F_p$, so the prime-to-$p$ part of $J_1(13)(K)_{tors}$ injects into $J_1(13)(\F_{p})$.

As $3$ and $47$ are inert in $\Q(\sqrt 5)$, $|J_1(13)(\F_9)|=3\cdot 19$ and $|J_1(13)(\F_{47^2})|=2^8\cdot 7^2 \cdot 19^2$, we conclude $J_1(13)(\Q(\sqrt 5))_{tors}=J_1(13)(\Q)_{tors}$.

As $3$ and $5$ are inert in $\Q(\sqrt {-7})$, $|J_1(13)(\F_9)|=3\cdot 19$ and $|J_1(13)(\F_{25})|= 19^2$, we conclude $J_1(13)(\Q(\sqrt{-7}))_{tors}=J_1(13)(\Q)_{tors}$.

As $3$ and $11$ are inert in $\Q(\sqrt {2})$, $|J_1(13)(\F_9)|=3\cdot 19$ and $|J_1(13)(\F_{121})|= 7^2\cdot 19^2$, we conclude $J_1(13)(\Q(\sqrt{2}))_{tors}=J_1(13)(\Q)_{tors}$.

As $5$ and $29$ are inert in $\Q(\sqrt {-2})$, $|J_1(13)(\F_{25})|=19^2$ and $|J_1(13)(\F_{29^2})|=2^6\cdot 3^2\cdot 19 \cdot 61$, we conclude $J_1(13)(\Q(\sqrt{-2}))_{tors}=J_1(13)(\Q)_{tors}$.

As $3$ splits and $41$ is inert in $\Q(\sqrt {13})$ and $\Q(\sqrt {-11})$, and $|J_1(13)(\F_3)|=19$ and $|J_1(13)(\F_{41^2})|=2^6\cdot 7^4\cdot 19$, we conclude $J_1(13)(\Q(\sqrt{-11}))_{tors}=J_1(13)(\Q(\sqrt{13}))_{tors}$ $=J_1(13)(\Q)_{tors}$.

As $5$ and $17$ are inert in $\Q(\sqrt {3})$, $|J_1(13)(\F_{25})|=19^2$ and $|J_1(13)(\F_{17^2})|=2^6\cdot 3^2\cdot 7 \cdot 19$, we conclude $J_1(13)(\Q(\sqrt{3}))_{tors}=J_1(13)(\Q)_{tors}$.

As $17$ splits and $41$ is inert in $\Q(\sqrt {-15})$, $|J_1(13)(\F_{17})|=2^2\cdot 3 \cdot 19$ and $|J_1(13)(\F_{41^2})|=2^6\cdot 7^4\cdot 19$, we conclude $J_1(13)(\Q(\sqrt{-15}))_{tors}=J_1(13)(\Q)_{tors}$.

Note that Reichert \cite{re} already found an elliptic curve with torsion $\Z/13 \Z$ over $\Q(\sqrt{17})$, $$y^2 = x^3 -( 4323+1048\sqrt{17})x +227630+55208\sqrt{17},$$ where $(-49-12\sqrt{-7},- 296 -72\sqrt{-7})$ is a point of order $13$.

\end{proof} The following theorem solves the problem for $\Z/14 \Z$. \newtheorem{tm4}[tm]{Theorem} \begin{tm4} \label{t4} The quadratic field $K$ with smallest $|\Delta |$ such that $\Z/14 \Z$ appears as a torsion group over $K$ is $K=\Q (\sqrt{-7})$. \end{tm4} \begin{proof}

Note that from \cite{Naj1} we can see that $\Z/14 \Z$ does not appear over $\Q(i)$ and $\Q(\sqrt{-3})$. $X_1(14)(\Q)_{tors}\simeq \Z/6 \Z$, and all the points in $X_1(14)(\Q)_{tors}$ are cusps. From \cite[Lemme 2.2, p. 32]{Rab}, we see that $X_1(14)(\Q(\sqrt{-7}))_{tors}\simeq \Z/2 \Z \oplus \Z/6 \Z$, and $X_1(14)(K)_{tors}\simeq \Z/6 \Z$ for all other quadratic fields $K$. We compute $rank(X_1(14)(\Q(\sqrt{5})))=0$, so all the points on $X_1(14)(\Q(\sqrt{5}))$ are cusps.

By examining torsion points on $X_1(14)(\Q(\sqrt{-7}))$ that are not rational, we find noncuspidal points on $X_1(14)(\Q(\sqrt{-7}))$. This induces the curve $$y^2 +\frac{63+\sqrt{-7}}{56}xy+\frac{11+\sqrt{-7}}{112}y=x^3+\frac{11+\sqrt{-7}}{112}x^2,$$ where the point $(0,0)$ is of order $14$.

\end{proof}

The following theorem solves the problem for $\Z/15 \Z$.

\newtheorem{tm5}[tm]{Theorem} \begin{tm5} \label{t5} The quadratic field $K$ with smallest $|\Delta |$ such that $\Z/15 \Z$ appears as a torsion group over $K$ is $K=\Q (\sqrt{5})$. \end{tm5} \begin{proof} Note that from \cite{Naj1} we can see that $\Z/15 \Z$ does not appear over $\Q(i)$ and $\Q(\sqrt{-3})$. $X_1(15)(\Q)_{tors}\simeq \Z/4 \Z$, and all the points in $X_1(15)(\Q)_{tors}$ are cusps. The only quadratic fields $K$ such that $X_1(15)(K)_{tors}\not\simeq \Z/4 \Z$ are $K=\Q(\sqrt{-3}),\ \Q(\sqrt{5})$ and $\Q(\sqrt{-15})$. One easily checks that the points on $X_1(15)(\Q(\sqrt{-3}))_{tors}$ are cusps, while points on $X_1(15)(\Q(\sqrt{5}))_{tors}$ and $X_1(15)(\Q(\sqrt{-15}))_{tors}$ are not. From a noncuspidal point on $X_1(15)(\Q(\sqrt{5}))_{tors}$ we obtain an elliptic curve (already found by Reichert \cite{re}) $$y^2 = x^3 + (281880\sqrt{5} - 630315)x + 328392630-146861640\sqrt{5},$$ with a point $(264\sqrt{5} - 585, 5076\sqrt{5} - 11340)$ of order $15$. \end{proof} \newtheorem{rem}[tm]{Remark} \begin{rem} Note that the two example curves in Theorems \ref{t4} and \ref{t5} prove that \cite[Example 2.5]{km} is wrong. We believe that the reason for this is that the bad reduction of the modular curve $X_0(N)$ (for $N=14, 15$) at the primes $7$ and $5$ respectively was not taken into account. \end{rem}

We turn to the group $\Z/16  \Z$.

\newtheorem{tm6}[tm]{Theorem} \begin{tm6} The quadratic field $K$ with smallest $|\Delta |$ such that $\Z/16  \Z$ appears as a torsion group over $K$ is $K=\Q (\sqrt{-15})$. \end{tm6} \begin{proof} Note that from \cite{Naj1} we can see that $\Z/16  \Z$ does not appear over $\Q(i)$ and $\Q(\sqrt{-3})$. The torsion of the Jacobian of $X_1(16)(\Q)$, $J_1(16)(\Q)_{tors}$ is isomorphic to $\Z/2 \Z \oplus \Z/10 \Z$. All of these torsion points are induced by cusps of $X_1(16)$. We compute that $rank(J_1(16)(\Q(\sqrt d)))=0$ for $d=5,-7,2,-2,-11,3$ and $13$, and in a similar manner like in Theorem \ref{t3}, we obtain $J_1(16)(\Q(\sqrt d))_{tors}=J_1(16)(\Q)_{tors}$ for all the mentioned values $d$, with the exception of $d=2$. For $d=2$, we obtain $J_1(16)(\Q)_{tors}=\Z/2 \Z \oplus \Z/2 \Z \oplus \Z/10 \Z$, but again all the points on the curve $X_1(16)(\Q(\sqrt 2))$ are cusps.

The elliptic curve (taken from \cite{re}) $$y^2=x^3+272133x+41173974,$$ has a point $(3-144\sqrt{-15}, -6480-432 \sqrt{-15})$ of order $16$ over $\Q(\sqrt{-15})$. \end{proof}

\newtheorem{tm7}[tm]{Theorem} \begin{tm7} The quadratic field $K$ with smallest $|\Delta |$ such that $\Z/18 \Z$ appears as a torsion group over $K$ is $K=\Q (\sqrt{33})$. \end{tm7} \begin{proof} Note that from \cite{Naj1} and \cite{Naj2} we can see that $\Z/18 \Z$ does not appear over $\Q(i)$ and $\Q(\sqrt{-3})$.  We can immediately disregard the fields $\Q(\sqrt{d})$ for $d=2,-19,-31,17,21$ as a consequence of \cite[Proposition 2.4 (i)]{km}, $d=-11,13,21$ as a consequence of \cite[Proposition 2.4 (ii)]{km} and $d=-7,-15$ as a consequence of \cite[Proposition 2.4 (iii)]{km}. One is left to deal with the cases $d=-2,3,23,6,-6$. Computing $rank(J_1(18)(\Q(\sqrt d)))=0$  and $J_1(18)(\Q(\sqrt d))_{tors}=J_1(18)(\Q(\sqrt d))=\Z_{21}$ for all these cases, we get that there are no noncuspidal points on $X_1(18)$ for any of these fields.

For an example of an elliptic curve with torsion $\Z/18 \Z$ over $\Q(\sqrt{33})$ we take another example of Reichert \cite{re}, $$y^2 = x^3 + (28296\sqrt{33} - 162675)x + 35441118-6168312\sqrt{33},$$ where $(147-24\sqrt{33},540 -108\sqrt{33})$ is a point of order $18$.

\end{proof}

\newtheorem{tm8}[tm]{Theorem} \begin{tm8} \label{t8} The quadratic field $K$ with smallest $|\Delta |$ such that $\Z/2 \Z \oplus \Z/10 \Z$ appears as a torsion group over $K$ is $K=\Q (\sqrt{-2})$. \end{tm8} \begin{proof} Note that from \cite{Naj1} we can see that $\Z/2 \Z \oplus \Z/10 \Z$ does not appear over $\Q(i)$ and $\Q(\sqrt{-3})$. All the points on $X_1(2,10)(\Q(\sqrt{-7}))=X_1(2,10)(\Q) \simeq \Z/6 \Z$ are cusps, as are  all the points on $X_1(2,10)(\Q(\sqrt{5}))\simeq \Z/2 \Z \oplus \Z/6 \Z$. Both $X_1(2,10)(\Q(\sqrt{-7}))$ and $X_1(2,10)(\Q(\sqrt{5}))$ have rank 0, while $X_1(2,10)$ $(\Q(\sqrt{-2}))$ has rank 1.

For an elliptic curve with torsion $\Z/2 \Z \oplus \Z/10 \Z$ over $\Q(\sqrt{-2})$ we take a nontorsion point on $X_1(2,10)(\Q(\sqrt{-2}))$ and obtain the curve $$y^2 + \frac 5 {11}xy + \frac 6 {121}y = x^3 + \frac 6 {121}x^2,$$ where the points $((-2-8\sqrt{-2})/121,(20\sqrt{-2} - 28)/1331)$ and $(6/11,-72/121)$ generate the torsion group $\Z/2 \Z \oplus \Z/10 \Z$.

\end{proof}

\newtheorem{tm11}[tm]{Theorem} \begin{tm11} \label{t11}

The quadratic field $K$ with smallest $|\Delta |$ such that $\Z/2 \Z \oplus \Z/12 \Z$ appears as a torsion group over $K$ is $K=\Q (\sqrt{13})$. \end{tm11} \begin{proof} Note that from \cite{Naj1} we can see that $\Z/2 \Z \oplus \Z/12 \Z$ does not appear over $\Q(i)$ and $\Q(\sqrt{-3})$. The field with the smallest $|\Delta |$ such that $X_1(2,12)$ has positive rank over it is $\Q(\sqrt{13})$. Note that for all quadratic fields $K$ except for $K=\Q(i)$, $\Q(\sqrt{3})$ or $\Q(\sqrt{-3})$, $X_1(2,12)(K)_{tors}=X_1(2,12)(\Q)_{tors}\simeq \Z/4 \Z$ holds (see \cite[Lemme 2.4 p.40]{Rab}). For $K=\Q(i)$, $\Q(\sqrt{3})$ and $\Q(\sqrt{-3})$ all the points on $X_1(2,12)(K)_{tors}$ are cusps.

For an elliptic curve with torsion $\Z/2 \Z \oplus \Z/12 \Z$ over $\Q(\sqrt{13})$ we take a nontorsion point on $X_1(2,12)(\Q(\sqrt{13}))$ and obtain the curve $$y^2 + ( 134025-37172\sqrt{13})xy + (47915630355840-13289404780320\sqrt{13} )y=$$ $$ x^3 + (3775925760\sqrt{13} - 13614293940)x^2,$$ where the points $$(954691712-264783840\sqrt{13},42132429627392\sqrt{13}-151910635381440)$$ and $$(3993089880-1107483870\sqrt{13}, 176222937989280\sqrt{13} - 635380838833260)$$ generate the torsion group $\Z/2 \Z \oplus \Z/12 \Z$. \end{proof}

\section{The Method of Mazur for $X_1(13)$ and $X_1(18)$}

As mentioned earlier, the curves $X_1(13)$ and $X_1(18)$ are the hardest to deal with. One can expect that for a large number of fields, $2$-descent will not be able to prove the finiteness of $J_1(N)$. We use a different method to obtain a criterion for the finiteness of $J_1(N)$ that will be satisfied for many imaginary quadratic fields.

In \cite{Kam2} the first author carried out a 19-Einstein prime descent on $J_1 (13)$ over a quadratic imaginary field $K$. Here we perform an analogous task, a 7-Eisenstain prime descent on $J_1 (18)$ over a quadratic imaginary field. We describe the results of these descents in Theorem \ref{kt1}.

When $N=13$ or 18 the curve $X_1 (N)$ is of genus 2. We let $G$ be the Galois group of the cover $X_1(N) \to  X_0(N)$, so $G \approx ({\Bbb Z} /N {\Bbb Z})^*/(\pm 1)$.  In each case the Hecke algebra ${\Bbb T}$ (see \cite{Kam2}) is isomorphic to ${\Bbb Z}[G] \approx  {\Bbb Z}[\zeta_3]$. When $N=13$ the ${\Bbb Q}(\zeta_N)^+$-rational cuspidal group $C$ has order 19, and when $N=18\; C$ has order 21. The Eisenstein ideal $I$ is the ideal of ${\Bbb T}$ that annihilates $C$. We set $q=19$ if $N=13$, and $q=7$ if $N=18$. The $q$-Eisenstein prime $\pi$ is the prime ideal of ${\Bbb T}$, above $q$, and in the support of $I$. It is the annihilator of the $q$-Sylow subgroup $C_q$ of $C$. We write $\varepsilon$ for the character via which $G$ acts on $C_q$, and we write ${\Bbb Q}_\varepsilon$ for the field cut out by $\varepsilon$ (by identifying $G$ with Gal$({\Bbb Q}(\zeta_N)^+/{\Bbb Q}))$.

Let $K$ be an imaginary quadratic field with associated character $\psi$, and suppose $\psi$ is disjoint from ${\Bbb Q}_\varepsilon$.  We write $E$ for the composition $K \cdot {\Bbb Q}_\varepsilon$. We decompose the $q$-Sylow subgroup ${\bf C}$ of the ideal class  group of $E$ into eigenspaces under the action of Gal$(E/{\Bbb Q})$, $$ {\bf C} = \oplus \: {\bf C}(\varepsilon^a \cdot \psi^b) $$ where ${\bf C} (\varepsilon^a \psi^b) = \{c \in {\bf C} : \sigma c = \varepsilon^a \psi^b(\sigma) \cdot c\;  \forall \sigma \in$ Gal$(E/Q)\}$. We note that ${\bf C}(\varepsilon^{-a} \psi^{-b})\neq0$ if and only if ${\Bbb B}_{1, \varepsilon^{a} \psi^{b}} \equiv 0 \; (\mbox{mod }\pi)$, where ${\Bbb B}_{1,\chi}$ is the first generalized Bernoulli number.

\newtheorem{kt1}[tm]{Theorem} \begin{kt1} \label{kt1} If $N=18$ we assume that 2 doesn't split in $K$. If ${\bf C}(\psi)$ and ${\bf C}(\psi \varepsilon^{-1})$ are both zero then the Mordell-Weil group $J_1(N)(K)$ is finite. \end{kt1}

\begin{proof} When $N=13$ this is Theorem 6.1 of \cite{Kam2}. For the remainder of the proof we will assume that $N=18$. The Hecke algebra is ${\Bbb Z}[\zeta_3]$, so the ideal (7) splits into a product $(7)= \pi \cdot \overline{\pi}$ of 2 primes. One of these is the Eisenstein prime $\pi$, and the other is the annihilator $\overline{\pi}$ of the ${\Bbb Q}$-rational cuspidal subgroup of order 7. The 7-torsion in $J_1(N)$ has a corresponding decomposition $J[7]=J[\pi] \oplus J[\overline{\pi}]$ into the direct sum of the $\pi$ and $\overline{\pi}$-torsion subgroups. As a Gal$(\overline{\Bbb Q}/{\Bbb Q})$-module $C$ is isomorphic to ${\Bbb Z}/7 {\Bbb Z}[\varepsilon]$, the twist of the constant Galois module by the character $\varepsilon$. By the Eichler-Shimura relation the quotient $J[\pi]/C$ is isomorphic to  ${\boldsymbol\mu}_7$, the Gal$(\overline{\Bbb Q}/{\Bbb Q})$-module of 7th roots of unity (see \cite{Kam2} for details).

We will carry out a $\pi$-descent in the style of Mazur. One might think that the descent would be made easier if we utilize the potential good reduction at 3, and then take Gal$({\Bbb Q}_\varepsilon/{\Bbb Q})$ - invariants (which is easy since $[{\Bbb Q}_\varepsilon : {\Bbb Q}]$ is prime to 7). However, it turns out that one doesn't actually gain anything by doing this. As we shall see the contribution to the relevant cohomology groups coming from the prime 3 turns out to be trivial.

Since ${\Bbb T}={\Bbb Z}[\zeta_3]$ is a P.I.D., the ideal $\pi$ is principal, say $\pi =(\eta)$.  Let $\triangle_K = \{\wp : \wp|2$ or 3 in ${\cal O}_K\}$, and let $S=$ Spec ${\cal O}_K - \triangle_K$. Write $A$ for the open subscheme of the N\'{e}ron model ${\cal A}_{/{\mbox{Spec }} {\cal O}_{K}}$ whose fibers are connected. We examine the sequence \begin{equation} \label{ke1} 0 \to A[\eta] \to A \stackrel{\cdot \eta}{\to} A \to 0 \ldots \end{equation} over the base Spec ${\cal O}_K$. The group scheme $A[\eta]$ is a quasi-finite flat group scheme over Spec ${\cal O}_K$. We will also need to study sequence (\ref{ke1}) over the base $S$ where $A$ is an abeliam scheme, and $A[\eta]$ is a finite flat group scheme.

\newtheorem{kt2}[tm]{Lemma} \begin{kt2} \label{kt2} (a) There is a short exact sequence of finite flat group schemes $$ 0 \to C \to A[\eta] \to {\boldsymbol \mu}_7 \to 1. $$ over the base $S$.\\ (b)	There is a short exact sequence of quasi-finite flat group schemes $$ 0 \to \overline{C} \to A[\eta] \to \overline{{\boldsymbol \mu}_7} \to 1 $$ over the base Spec ${\cal O}_K$. Here $\overline{C}$ and $\overline{{\boldsymbol\mu}_7}$ are extensions of $C_{/S}$ and ${\boldsymbol \mu}_{7/S}$ to Spec ${\cal O}_K$. \end{kt2}

\begin{proof} Part (a) follows from the comments at the beginning of the proof of Theorem \ref{kt1}, together with the Oort-Tate classification of group schemes of prime order.\\ Part (b) is just obtained by taking the closure of the group schemes in (a). \end{proof} We use sequence \ref{kt2}(b) to analyze the $f.p.q.f$. cohomology of sequence (\ref{ke1}), i.e., we study the exact sequence $$A(Spec\ {\cal O}_K) \stackrel{\cdot \eta}{\to} A (Spec\ {\cal O}_K) \to  H^1(Spec\ {\cal O}_K, A[\eta])$$ with the goal of showing that the right hand group is trivial. This implies that $A$(Spec ${\cal O}_K)$ is finite. Together with the fact that $A$(Spec ${\cal O}_K)$ is of finite index in ${\cal A}$(Spec ${\cal O}_K) = J_1(18)(K)$ we see that the latter group is finite.

To compute $H^1$(Spec ${\cal O}_K, A[\eta])$ we compute $H^1$(Spec ${\cal O}_K, \overline{C})$ and $H^1$(Spec ${\cal O}_K, \overline{{\boldsymbol \mu}_7})$. In order to analyze the latter we use the sequence $$ 1 \to \overline{{\boldsymbol \mu}_7} \to {\boldsymbol \mu}_7 \to {\boldsymbol \mu}_{7^{\cdot}} \to 1 $$ where ${\boldsymbol \mu}_{7^{\cdot}}$ is a skyscraper sheaf concentrated at the points of $\triangle$. As usual, $H^1$(Spec ${\cal O}_K, {\boldsymbol \mu}_7)$ fits into a Kummer sequence $$ 1 \to {\cal O}^*_K/({\cal O}^*_K)^7 \to H^1(\mbox{Spec }{\cal O}_K, {\boldsymbol \mu}_7) \to \mbox{Pic}({\cal O}_K)[7] \to 1. $$ The left hand group is trivial since $K$ is an imaginary quadratic field, and the right is trivial since the class number of $K$ is assumed to be relatively prime to 7 (i.e., ${\bf C}(\psi)=0)$. Thus $H^1$(Spec ${\cal O}_K, {\boldsymbol \mu}_7)=0$.\\ Now $H^\circ(Spec\ {\cal O}_K, {\boldsymbol \mu}_{7^{\cdot}}) ={{\oplus }\atop{\wp \in \triangle }} {\boldsymbol \mu}_7({\cal O}_K/\wp),$ and so this group is also trivial. Thus, the sequence $$ \to H^\circ({\mbox{Spec }}{\cal O}_K, {\boldsymbol \mu}_{7^{\cdot}}) \to H^1({\mbox{Spec }}{\cal O}_K, \overline{{\boldsymbol \mu}_7}) \to H^1({\mbox{Spec }}{\cal O}_K, {\boldsymbol \mu}_7) $$
 shows that the middle group is trivial.

We turn our attention to the group $H^1$(Spec ${\cal O}_K, \overline{C})$. Let $\triangle_E = \{\wp \in$ Spec ${\cal O}_E : \wp|2$ or $3\}$, and let $T=$Spec ${\cal O}_E-\triangle_E$. Mazur [6] shows that $H'$(Spec ${\cal O}_K, \overline{C})$ injects into $H^1(S, C)$, and the latter group injects into $H'(T,C_{/T})^{{\mbox{Gal}}(T/S)} \approx H^1(T, {\Bbb Z}/7 {\Bbb Z})^{{\mbox{Gal}}(T/S)}$ (since $C_{/T} \approx {\Bbb Z}/7 {\Bbb Z}$). The latter injection follows from the Hochschild-Serre spectral sequence $H^p$(Gal($T/S),$ $H^q(T, C_{/T})) \Rightarrow H^{p+q} (S, C)$ which degenerates since Gal$(T,S)$ has order prime to 7.

We examine the long exact relative cohomology sequence \begin{equation} \label{ke2} \begin{array}{ccc} 0 \to H^1({\mbox{Spec }}{\cal O}_E, {\Bbb Z}/7 {\Bbb Z}) \to H^1(T, {\Bbb Z}/7 {\Bbb Z}) \to H^2_{\triangle}({\mbox{Spec }} {\cal O}_E, {\Bbb Z}/7 {\Bbb Z})\\ \to H^2({\mbox{Spec }}{\cal O}_E, {\Bbb Z}/7 {\Bbb Z}) \ldots \end{array} \end{equation} Now $H^2_\triangle$ (Spec ${\cal O}_E, {\Bbb Z} /7 {\Bbb Z}) \approx {{\oplus}\atop{\wp \in \triangle}} H^2\cdot$(Spec ${\cal O}_{E, \wp}, {\Bbb Z}/{\Bbb Z})$, which is isomorphic to the Pontryagin dual ${{\oplus}\atop{\wp \in \triangle}} H^1($Spec ${\cal O}_{E, \wp}, {\boldsymbol \mu}_7)^*$ by local flat duality.

The map $H^2_\triangle$ (Spec ${\cal O}_{E, \wp}, {\Bbb Z}/7{\Bbb Z}) \to H^2/$(Spec ${\cal O}_E, {\Bbb Z} /7 {\Bbb Z})$ of sequence (\ref{ke2}) is dual to the map $H^1$(Spec ${\cal O}_E, {\boldsymbol \mu}_7) \stackrel{\alpha}{\to} H^1$(Spec ${\cal O}_{E, \wp}, {\boldsymbol \mu}_7)$, where $\wp$ is the unique prime of $E$ dividing 2. Under $\alpha$ the group $({\cal O}^*_E)/({\cal O}^*_E)^7$ maps surjectively onto $k^*_\wp/(k^*_\wp)^7$ (where $k_\wp = {\cal O}_E/\wp)$. Then $H^2_\triangle$ (Spec ${\cal O}_E, {\Bbb Z} /7 {\Bbb Z})$ maps injectively to $H^2$(Spec ${\cal O}, {\Bbb Z}/7 {\Bbb Z})$, which tells us (from (\ref{ke2})) that $H^1$(Spec ${\cal O}_E, {\Bbb Z}/7 {\Bbb Z})$,  surjects onto $H^1(T, {\Bbb Z}/7 {\Bbb Z})$, i.e. $$ \begin{array}{lll} H^1(\mbox{Spec } {\cal O}_E, {\Bbb Z} /7 {\Bbb Z}) \approx H^1 (T, {\Bbb Z} /7 {\Bbb Z})\\ \approx|\\ \mbox{Hom }({\cal C}l_E, {\Bbb Z} /7 {\Bbb Z}). \end{array} $$ Taking Gal$(T/S)$-invariants and using the assumption ${\cal C}(\psi \varepsilon^{-1}) = 0$ we find $H^1(T, {\Bbb Z} /7 {\Bbb Z})^{\mbox{Gal }(T/S)} = 0$, as desired. This completes the proof that $J_1 (N)(K)$ is finite.

\end{proof}

Finally, we wish to use the results of the descents to study $K$-rational points on $X_1(N)$. Let $R =$ Spec ${\Bbb Z}$ if $N=13$, or Spec ${\Bbb Z}[1/3]$ if $N=18$. Mazur's argument (\cite{Maz2}, 3.1) works mutatis mutandis in this case and proves the following.

\vskip 1pc

\newtheorem{kt3}[tm]{Proposition} \begin{kt3} \label{kt3} The morphism $f: X_1(N)^{\mbox{\tiny {smooth}}}_{/R} \to J_1(N)$ is a formal immersion away from characteristic 2 and 13 if $N=13$, and away from characteristic 2 if $N=18$. \end{kt3}

We recall that $K$ is a quadratic imaginary field such that $J_1(N)(K)$ is finite. If $N=13$ let $p=3$ or 5, and if $N=18$ let $p=5$ or 7, and assume that $p$ splits or ramifies in $K$. We let $\wp$ be a prime above $p$ in ${\cal O}_K$, and write $k$ for ${\cal O}_K/_{\wp}$. Proposition \ref{kt3}, together with the argument of \cite{Maz2}, Corollary 4.3, shows that if $E$ is an elliptic curve over $K$ with a $K$-rational point $P$ of order $N$ then $E$ has potentially good reduction at $\wp$. The reduction cannot be good because $N$ is too large for the reduction of $P$ to exist on $E_{/k}$.  The reduction must therefore be additive.  However, $N$ is also too large to divide $[E : E^\circ]$ in the additive case, which forces $P$ to specialize to $E^\circ_{/k} \approx G_a$. This is clearly impossible since gcd$(N,p)=1$, so the pair $(E, P)$ cannot exist.

\section{Interplay of rank and torsion}

In this section we will study how often a torsion group appears over a fixed quadratic field if it appears at all. All three possibilities will happen: some groups will appear infinitely often over a fixed quadratic field $K$, some finitely, and some can appear finitely and infinitely, depending on $K$.

\newtheorem{tm9}[tm]{Theorem} \begin{tm9} \label{t9} Suppose that over a quadratic field $K$, $\Z/13 \Z, \Z/16  \Z$ or $\Z/18 \Z$ appears as torsion group of an elliptic curve. Then there are finitely many elliptic curves (up to isomorphism) over $K$ with that torsion group. \end{tm9} \begin{proof} As $X_1(N)$ is a curve of genus $>1$ for $N=13,16,18$, by Faltings' theorem, $X_1(N)(K)$ has finitely many points. \end{proof}

\newtheorem{tm10}[tm]{Theorem} \begin{tm10} Suppose that over a quadratic field $K$, $\Z/11 \Z, \Z/2 \Z\oplus \Z/10 \Z$ or $\Z/2 \Z\oplus \Z/12 \Z$ appears as torsion group of an elliptic curve. Then there are infinitely many elliptic curves over $K$ with that torsion group. \end{tm10} \begin{proof} Since $X_1(m,n)$ for $(m,n)=(1,11),(2,10),(2,12)$ are elliptic curves, one has to prove that there does not exist a field where $X_1(m,n)$ has a noncuspidal torsion point and rank 0.

As was mentioned in Theorem \ref{t2}, $X_1(11)(K)_{tors}=X_1(11)(\Q)_{tors}\simeq \Z/5 \Z$ and all the torsion points are cusps. Thus if $X_1(11)(K)$ has a noncuspidal point, it is of infinite order.

As was mentioned in Theorems \ref{t8} an \ref{t11}, all the points on $X_1(2,10)(K)_{tors}$ and $X_1(2,12)(K)_{tors}$ are cusps, for all quadratic fields $K$.

\end{proof} \newtheorem{tm12}[tm]{Theorem} \begin{tm12} \label{t12} (i) There are two elliptic curves with torsion $\Z/14 \Z$ over $\Q(\sqrt{-7})$. Over all other quadratic fields, if there exists one elliptic curve with torsion $\Z/14 \Z$, there exists infinitely many.\\ (ii) There is one elliptic curve with torsion $\Z/15 \Z$ over the fields $\Q(\sqrt{5})$ and $\Q(\sqrt{-15})$. Over all other quadratic fields, if there exists one elliptic curve with torsion $\Z/15 \Z$, there exists infinitely many. \end{tm12} \begin{proof} (i) As mentioned in Theorem \ref{t4}, $\Q(\sqrt{-7})$ is the only field over which the elliptic curve $X_1(14)$ has torsion that is larger than $X_1(14)(\Q)_{tors}$ and there are points on $X_1(14)(\Q(\sqrt{-7}))_{tors}$ that are not cusps. Moreover, $rank (X_1(14)(\Q(\sqrt{-7})))=0$, so there are finitely many points on $X_1(14)$ over $(\Q(\sqrt{-7}))$. By checking the curves generated by the torsion points of $X_1(14)(\Q(\sqrt{-7}))$ one easily sees that there are, up to isomorphism, exactly two elliptic curves with torsion $\Z/14 \Z$ over $\Q(\sqrt{-7})$. The curves have the following equations $$ y^2 + \frac{ 14+3\sqrt{-7}}{7}xy + \frac {-3+\sqrt{-7}}{7}y = x^3 + \frac {-3+\sqrt{-7}}{7}x^2,$$ and $$y^2 + \frac{ 7+2\sqrt{-7}}{7}xy + \frac{ 1+\sqrt{-7}}{7}y = x^3 + \frac{ 1+\sqrt{-7}}{7}x^2,$$ where $(0,0)$ is a point of order $14$ on both curves.

For all other quadratic fields, a noncuspidal point on $X_1(14)(\Q(\sqrt{-7}))$ will have infinite order.\\ (ii) As mentioned in Theorem \ref{t5}, the only fields over which the elliptic curve $X_1(15)$ has torsion larger than $X_1(15)(\Q)_{tors}$ are the fields $\Q(\sqrt{-3}), \Q(\sqrt{5})$ and $\Q(\sqrt{15})$. One checks that the points on $X_1(15)(\Q(\sqrt{-3}))_{tors}$ are cusps, while the points on $X_1(15)(\Q(\sqrt{5}))_{tors}$ and $X_1(15)(\Q(\sqrt{-15}))_{tors}$ are not. One computes that the $rank( X_1(15)(\Q(\sqrt{5})) )= rank ( X_1(15)(\Q(\sqrt{-15})) )=0$, proving that there are only finitely many curves with torsion $\Z/15\Z$ over these fields. By checking the curves generated by the torsion points of  $X_1(15)(\Q(\sqrt{5}))$ and $X_1(15)(\Q(\sqrt{-15}))$, we find that over both of these fields there is exactly one curve with 15-torsion. Over $\Q(\sqrt{-15})$ this is the elliptic curve $$y^2 + \frac{145+7\sqrt{-15}}{128}xy + \frac{ 265+79\sqrt{-15}}{4096}y = x^3 + \frac{265+79\sqrt{-15}}{4096}x^2,$$ where $((95-9\sqrt{-15})/512, (255+65 \sqrt{-15})/16384)$ is a point of order $15$.

The only curve with 15-torsion over $\Q(\sqrt{5})$ can be found in the proof of Theorem \ref{t5}.

\end{proof}

One can examine not just the torsion group, but the whole Mordell-Weil group of these exceptional curves. One can easily compute that the rank of all four curves from Theorem \ref{t12} is equal to zero.

This proves the following corollary.

\newtheorem{tm13}[tm]{Corollary} \begin{tm13}

a) All elliptic curves with torsion $\Z/14 \Z$ over $\Q(\sqrt{-7})$ have rank zero.\\ b) All elliptic curves with torsion $\Z/15 \Z$ over $\Q(\sqrt{5})$ and $\Q(\sqrt{-15})$ have rank zero. \end{tm13}

Note that this is in stark contrast to what happens over the rationals; it is a widely believed conjecture that an elliptic curve with prescribed torsion can have arbitrary large rank. Moreover, for any torsion group appearing over the rationals, there exists an elliptic curve with that torsion group and rank at least $3$ (see \cite{Duj}).

An upper limit on the rank of an elliptic curve with given torsion over a fixed quadratic field can exist and be positive. For example, by Theorem \ref{t9}, there are only finitely elliptic curves with torsion $\Z/13 \Z$ over $\Q(\sqrt{193})$, and by \cite[Theoreme 2.4]{Rab}, there exists a curve with rank at least 2.

\section{Number of fields having given torsion}

It is natural to wonder what is the density of the fields having some fixed group as a torsion group of elliptic curves. For torsion groups such that the curve $X_1$ inducing them is an elliptic curve, this will depend on the rank of the quadratic twists of $X_1$. If $T=\Z/m\Z \oplus \Z/nZ$ is a torsion group that is generated by an elliptic modular curve $X_1(m,n)$, then for all except finitely many quadratic fields $K$, $T$ will appear as a torsion group over $K$ if and only if $X_1$ has positive rank over $K$. This means that the standard conjecture \cite[Corollary 7.4]{rs}, implies that the group $T$ will appear over $1/2$ of the quadratic fields, when ordered by discriminant. However, the best that one can at the moment unconditionally prove is an asymptotical lower bound on the number of fields $\Q(\sqrt d),\ |d|\leq X$ for some bound $X$, having or not having $T$ as a torsion group.

\newtheorem{tm14}[tm]{Proposition} \begin{tm14} a) Let $T=\Z/14 \Z,\ \Z/15 \Z,\ \Z/2 \Z \oplus \Z/10 \Z$ or $\Z/2 \Z \oplus \Z/12 \Z$. Then the number of quadratic fields $\Q(\sqrt d)$ such that $0<|d|\leq X$ and $T$ does not appear as a torsion group of an elliptic curve over $\Q(\sqrt d)$ is $\gg X/\log X$.\\ b) The number of quadratic fields $\Q(\sqrt d)$ such that $0<|d|\leq X$ and $\Z/11 \Z$ does not appear as a torsion group of an elliptic curve over $\Q(\sqrt d)$ is $\gg X/(\log X)^{1-\alpha}$ for some $0<\alpha<1$.\\ c) Let $\epsilon>0$ and $T=\Z/11 \Z,\ \Z/14 \Z,\ \Z/15 \Z,\ \Z/2 \Z \oplus \Z/10 \Z$ or $\Z/2 \Z \oplus \Z/12 \Z$. Then the number of quadratic fields $\Q(\sqrt d)$ such that $0<|d|\leq X$ and $T$ appears as a torsion group (for infinitely many non-isomorphic curves) over $\Q(\sqrt d)$ is $\gg X^{1-\epsilon}$.\\ \end{tm14} \begin{proof} a) As $X_1(14),\ X_1(15),\ X_1(2,10)$ and $X_1(2,12)$ are elliptic curves, by \cite[Corollary 3]{Ons} we see that $\gg X/\log X$ of the twists of the appropriate curve $X_1$ will have rank $0$, and thus there will be no elliptic curves with given torsion over the corresponding quadratic field.

b) The proof is the same as in a) with the difference that since $X_1(11)$ does not have a rational $2$-torsion point, we can apply the stronger result \cite[Corollary 3]{Ono} and thus prove the theorem.

c) By \cite[Theorem 1]{pp}, $\gg X^{1-\epsilon}$ of the twists of each of the elliptic curves $X_1$ will have positive rank, and thus there will exist infinitely many elliptic curves with the appropriate torsion over this field.

\end{proof}

If $T$ is a torsion group induced by a hyperelliptic curve, then one expects this torsion group to not appear much more often than it appears. Indeed, one can prove that for a large density of quadratic fields, if the fields $\Q(\sqrt d)$ are ordered by the largest prime appearing in $d$, $\Z/18 \Z$ will not appear as a torsion group.

We will use the following Proposition of Kenku and Momose \cite[Proposition (2.4)]{km}.

\newtheorem{tm15}[tm]{Proposition} \begin{tm15} \label{t15} Let $K$ be a quadratic field. If the field satisfies one of the conditions (i), (ii) or (iii) listed below, then there are no elliptic curves over $K$ with torsion $\Z/18 \Z$.\\ (i) $3$ remains prime in $K$.\\ (ii) $3$ splits in $K$ and $2$ does not split in $K$.\\ (iii) $5$ or $7$ ramify in $K$. \end{tm15} We now define the function $\psi$ that will give us an ordering of quadratic fields that will be more suitable for our purposes. \newtheorem{tm16}[tm]{Definition} \begin{tm16} Let $d$ be a square-free integers and write $d_1=(-1)^{\alpha_0}2^{\alpha_1}\ldots p_k^{\alpha_k}$. We define the bijection $\psi$ from the set of all quadratic fields to the positive integers, by $\psi(\Q(\sqrt{d}))=(\alpha_k\ldots \alpha_0)_2$. \end{tm16}

We now prove the following.

\newtheorem{tm17}[tm]{Theorem} \begin{tm17} Define $N_t=\{\Q(\sqrt{d})|\ \Z/18 \Z$ does not appear as a torsion group of an elliptic curve over $\Q(\sqrt{d})$ and $\psi(\Q(\sqrt{d}))\leq t\}$ and $A_t=\{\Q(\sqrt{d})|$ $\psi(\Q(\sqrt{d}))\leq t\}$. Then $$\lim_{t\rightarrow \infty}\frac{N_t}{A_t}\geq \frac{55}{64}.$$ \end{tm17} \begin{proof} When ordered by the $\psi$ function, $1/4$ of quadratic fields will satisfy Proposition \ref{t15} (i), and $3/16$ will satisfy (ii). Since the set of fields satisfying (i) and the set of fields satisfying (ii) are disjoint, $7/16$ of the fields will satisfy either one. The condition (iii) is satisfied by $3/4$ of the fields and the probability that this condition is satisfied does not depend on whether (i) or (ii) are satisfied, and thus we conclude that at least $7/16+3/4 -7/16\cdot 3/4=55/64$ of quadratic fields (when ordered by $\psi$) does not have $\mathbb Z/18 \Z$ as a torsion subgroup appearing over it. \end{proof}

\textbf{Acknowledgements.} We are grateful to Andrej Dujella and Matija Kazalicki for helpful comments. The second author was financed by the National Foundation for Science, Higher Education and Technological Development of the Republic of Croatia.

\noindent \small{SHELDON KAMIENNY}\\

\noindent \small{DEPARTMENT OF MATHEMATICS, UNIVERSITY OF SOUTHERN CALIFORNIA, 3620 S. VERMONT AVE., LOS ANGELES, CA 90089-2532, U.S.A.}\\ \emph{E-mail address:} kamienny@usc.edu\\

\noindent \small{FILIP NAJMAN}\\

\noindent \small{MATHEMATISCH INSTITUUT, P.O. BOX 9512, 2300 RA LEIDEN, THE NETHERLANDS}\\ \emph{E-mail address:} fnajman@math.leidenuniv.nl\\

AND\\

\noindent \small{DEPARTMENT OF MATHEMATICS, UNIVERSITY OF ZAGREB, BIJENI\v CKA CESTA 30, 10000 ZAGREB, CROATIA}\\ \emph{E-mail address:} fnajman@math.hr

\end{document}